# Recursively squeezable sets are squeezable

by Fredric D. Ancel

**Abstract.** In work by Freedman [F2] and Freedman-Quinn [FQ] on the topology of 4-manifolds, null decompositions whose non-singleton elements are, in the terminology of [MOR], *recursively starlike-equivalent sets* of *filtration length 1* arise and are shown to be shrinkable. The main result of [MOR] is a general theorem covering these types of decompositions. It establishes the shrinkability of null decompositions whose non-singleton elements are recursively starlike-equivalent sets whose filtration lengths have a uniform finite upper bound. That result is the inspiration for this article. Here it is shown that the hypothesis of a uniform finite upper bound on filtration lengths is unnecessary. In outline: notions of *squeezable* subsets and *squashable* subsets of a compact metric space are defined. It is observed that starlike-equivalent sets are squeezable, and that any null decomposition of a compact metric space whose non-singleton elements are squeezable is shrinkable. It is also proved that a set is squeezable if and only if it is squashable, and that every recursively squashable set is squashable. It follows that any null decomposition of a compact metric space whose non-singleton elements are recursively squeezable is shrinkable. The latter theorem has as a corollary the main result of [MOR] with the hypothesis of a uniform finite upper bound on filtration lengths removed.

## 1. Introduction.

An *upper semi-continuous decomposition* of a topological space X is a partition $\mathcal{G}$ of X into compact subsets such that the associated quotient map X → X/$\mathcal{G}$ is a closed map. The quotient space X/$\mathcal{G}$ is called a *decomposition space.* A central question of decomposition space theory is: under what conditions is the associated quotient map X → X/$\mathcal{G}$ a near homeomorphism (i.e., approximable by homeomorphisms)? Decomposition space theory has played a key role in the effort to understand topological manifolds. A prime example is M. H. Freedman's proof in the 1980's of the 4-dimensional Poincaré Conjecture in the topological category [F1]. In that proof, a variant of traditional handlebody theory is performed in a 5-dimensional smooth proper h-cobordism using objects called *Casson handles* rather than topological open 2-handles. The proof succeeds in establishing that the h-cobordism is a product provided that Casson handles are homeomorphic to open 2-handles. It is then proved that such homeomorphisms exist. These homeomorphisms arise as approximations to quotient maps associated to decompositions. The homeomorphisms are outputs of limiting processes that fail in any obvious way to preserve the ambient differentiable structure, which explains why Freedman's result belongs to the topological category, and which perhaps sheds light on why the 4-dimensional Poincaré Conjecture in the differential category remains unresolved.

Subsequent work in [F2] and [FQ] extending Freedman's original theorems modified the sorts of decompositions under consideration and provided proofs that the associated quotient maps are near homeomorphisms. The new decompositions are of the type described in the abstract – null decompositions whose non-singleton elements are recursively starlike equivalent sets whose filtration lengths have a uniform finite upper bound. These are exactly the types of decompositions to which the main result of [MOR] applies. This article introduces the notions of *squeezable* and *squashable sets* and uses them to prove a theorem that generalizes the main result of [MOR] and requires no assumption of a finite upper bound on filtration lengths.

The author thanks Arunima Ray for a correspondence which led to improvements in the exposition.

## 2. Basic properties of decompositions.

**Notation.** Suppose X and Y are metrizable spaces with metrics $\rho$ and $\sigma$, respectively. For $S \subset X$ and $\epsilon > 0$, let $\mathcal{N}_\rho(S,\epsilon) = \{ x \in X : \rho(x,y) < \epsilon \text{ for some } y \in S \}$. Let C(X,Y) denote the set of all maps from X to Y, and let $\bar{\sigma}$ denote the supremum metric on C(X,Y) determined by $\sigma$. In particular, let C(X,X) denote the set of all maps from X to itself, and let $\bar{\rho}$ denote the supremum metric on C(X,X) determined by $\rho$. For $C \subset X$, let $\mathcal{H}_C(X)$ denote the set of all homeomorphisms $h : X \to X$ such that $h \mid C = id_C$, and let $cl(\mathcal{H}_C(X))$ denote the closure of $\mathcal{H}_C(X)$ in C(X,X). If $\pi : X \to Y$ is a map and $C \subset Y$ such that $\pi \mid \pi^{-1}(C) : \pi^{-1}(C) \to C$ is injective, we say that $\pi$ *injects over C.*

**Definition.** An *upper semi-continuous decomposition* of a topological space X is a partition $\mathcal{G}$ of X such that each element of $\mathcal{G}$ is compact and the quotient map $\pi : X \to X/\mathcal{G}$ is closed. (Here, $X/\mathcal{G}$ is given the *quotient topology:* a subset U of $X/\mathcal{G}$ is open if and only if $\pi^{-1}(U)$ is an open subset of X.) Note that the condition that the quotient map $\pi : X \to X/\mathcal{G}$ is closed is equivalent to the condition that the inverse of the quotient map satisfies the following continuity property: for every $y \in X/\mathcal{G}$ and every open neighborhood U of $\pi^{-1}(y)$ in X, there is an open neighborhood V of y in $X/\mathcal{G}$ such that $\pi^{-1}(V) \subset U$. From this point on, we will abbreviate the statement "$\mathcal{G}$ is an upper semi-continuous decomposition of X" to simply "$\mathcal{G}$ is a decomposition of X". For further information about decomposition spaces, we suggest the encyclopedic text [D].

**Remark.** If $\mathcal{G}$ is a decomposition of a (separable) metrizable space X, then $X/\mathcal{G}$ is also (separable) metrizable. This can be proved by showing that $X/\mathcal{G}$ satisfies the hypotheses of the relevant metrization theorem of Urysohn or Bing or Nagata-Smirnov. For further details, see Proposition 2 on page 13 of [D]. For the sake of simplicity, in this article we will focus on decompositions of compact metric spaces.

As noted in the introduction, a principal concern of decomposition space theory is determining whether the quotient map associated to a decomposition is a near homeomorphism. The main tool for identifying quotient maps that are near homeomorphisms is the Bing Shrinking Theorem. We now state the strong version of this result, first defining the relevant concepts.

**Definition.** Let $\pi : X \to Y$ be an onto map between compact metric spaces X and Y with metrics $\rho$ and $\sigma$, respectively, and let $\bar{\sigma}$ denote the supremum metric on $C(X,Y)$ determined by $\sigma$. We say that $\pi$ satisfies the *(strong) shrinking criterion* if:

for every $\delta > 0$ (and for every closed subset C of Y such that $\pi$ injects over C), there is a homeomorphism $h : X \to X$ such that $\text{diam}_\rho(h(\pi^{-1}(y))) < \delta$ for every $y \in Y$, $\bar{\sigma}(\pi, \pi \circ h^{-1}) < \delta$ (and $h \mid \pi^{-1}(C) = \text{id}_{\pi^{-1}(C)}$).

We say that $\pi$ is a *(strong) near homeomorphism* if for every $\epsilon > 0$ (and for every closed subset C of Y such that $\pi$ injects over C), there is a homeomorphism $f : X \to Y$ such that $\bar{\sigma}(\pi, f) < \epsilon$ (and $f \mid \pi^{-1}(C) = \pi \mid \pi^{-1}(C)$). Observe that since the limit in $C(X,Y)$ of onto maps must be onto, then each near homeomorphism is onto.

**The (Strong) Bing Shrinking Theorem.** An onto map $\pi : X \to Y$ between compact metric spaces is a (strong) near homeomorphism if and only if it satisfies the (strong) shrinking criterion.

**Remark.** R. H. Bing was the first to observe and exploit the fact that the shrinking criterion is a sufficient condition for $\pi$ to be a near homeomorphisms. However, he did not crystalize his observation into a formally stated theorem. This occurred later. (See page 42 of [D].) R. D. Edwards appears to have been the first to realize that the converse of Bing's observation is also true, and our statement of the (Strong) Bing Shrinking Theorem imitates the elegant formulation of the theorem given by Edwards in [E]. (We have added the parenthetical *strong* content.) Edwards sketched a brief and clever proof of the theorem in [E] that differed from previous proofs in that it uses the Baire Category Theorem in the space of maps from X to Y. (Proving the strong version of the theorem requires no additional effort.) An exposition of Edwards' proof can also be found on pages 23 to 25 of [D].

**Definition.** We call a decomposition $\mathcal{G}$ of a compact metric space X *(strongly) shrinkable* if the quotient map $\pi : X \to X/\mathcal{G}$ satisfies the (strong) shrinking criterion. Thus, the (Strong) Bing Shrinking Theorem implies the following result.

**Corollary.** If $\mathcal{G}$ is a decomposition $\mathcal{G}$ of a compact metric space X, then the quotient map $\pi : X \to X/\mathcal{G}$ spaces is a (strong) near homeomorphism if and only if $\mathcal{G}$ is (strongly) shrinkable. ∎

**Definition.** A collection $\mathcal{C}$ of subsets of a metric space X with metric $\rho$ is *null* if for every $\epsilon > 0$, $\{ C \in \mathcal{C} : \text{diam}_\rho(C) \geq \epsilon \}$ is finite. Thus, if $\mathcal{C}$ is null, then all but countably many elements of $\mathcal{C}$ are one-point sets. A decomposition of X which is a null collection is called a *null decomposition*. Note that if f : X → Y is an onto map between compact metric spaces and $\mathcal{C}$ is a null collection of subsets of X, then the uniform continuity of f implies that f($\mathcal{C}$) = { f(C) : C ∈ $\mathcal{C}$ } is a null collection of subsets of Y.

### 3. Squeezable sets

**Definition.** A subset A of a compact metric space X with metric ρ is a *squeezable* subset of X if A is compact and for every null decomposition $\mathcal{G}$ of X such that A ∈ $\mathcal{G}$, for every open neighborhood U of A in X, and for every δ > 0, there is a homeomorphism h ∈ $\mathcal{H}_{X-U}$(X) such that $\text{diam}_\rho$(h(A)) < δ and for every G ∈ $\mathcal{G}$, either $\text{diam}_\rho$(h(G)) < δ or h | G = $\text{id}_G$.

We state two theorems about squeezable sets. First: null decompositions whose non-singleton elements are squeezable sets are strongly shrinkable. The proof of this fact is a familiar argument to experts in decomposition space theory.

**Theorem 1.** If $\mathcal{G}$ is a null decomposition of a compact metric space X each of whose non-singleton elements is squeezable, then $\mathcal{G}$ is strongly shrinkable.

**Proof.** Let π : X → X/$\mathcal{G}$ denote the quotient map. Let ρ be a metric on X, let σ be a metric on X/$\mathcal{G}$, and let $\bar{\sigma}$ denote the supremum metric on C(X,X/$\mathcal{G}$) determined by σ. Let $G_1, G_2, G_3, \cdots$ be a list of the non-singleton elements of $\mathcal{G}$, and let $y_i$ = π($G_i$) for i ≥ 1. To prove that $\mathcal{G}$ is strongly shrinkable, let δ > 0 and let C be a closed subset of X/$\mathcal{G}$ that is disjoint from { $y_i$ : i ≥ 1 }. There is an n ≥ 1 such that $\text{diam}_\rho$($G_i$) < δ for all i > n. There are pairwise disjoint open subsets $V_1, V_2, \cdots, V_n$ of X/$\mathcal{G}$ − C such that $y_i$ ∈ $V_i$ and $\text{diam}_\sigma$($V_i$) < δ for 1 ≤ i ≤ n. For each i, 1 ≤ i ≤ n, since $G_i$ is squeezable, there is a homeomorphism $h_i$ ∈ $\mathcal{H}_{X-\pi^{-1}(V_i)}$(X) such that $\text{diam}_\rho$($h_i$($G_i$)) < δ and for every G ∈ $\mathcal{G}$, either $\text{diam}_\rho$($h_i$(G)) < δ or $h_i$ | G = $\text{id}_G$. Let U = $\cup_{1 \leq i \leq n} \pi^{-1}(V_i)$. Then a homeomorphism h : X → X is defined by setting h | $\pi^{-1}(V_i)$ = $h_i$ | $\pi^{-1}(V_i)$ for 1 ≤ i ≤ n and h | X − U = $\text{id}_{X-U}$. Clearly, $\text{diam}_\rho$(h(G)) < δ for each G ∈ $\mathcal{G}$. Also for 1 ≤ i ≤ n, and x ∈ $\pi^{-1}(V_i)$, since $h^{-1}(x)$ ∈ $\pi^{-1}(V_i)$, then σ(π(x),π∘$h^{-1}$(x)) ≤ $\text{diam}_\sigma$($V_i$) < δ; and for x ∈ X − U, σ(π(x),π∘$h^{-1}$(x)) = σ(π(x),π(x)) = 0. Thus, $\bar{\sigma}$(π,π∘$h^{-1}$) < δ. Finally, since $\pi^{-1}$(C) ⊂ X − U, then h | $\pi^{-1}$(C) = $\text{id}_{\pi^{-1}(C)}$. ∎

To state the second theorem about squeezable sets, we need some additional definitions.

**Definition.** A subset A of $\mathbb{R}^n$ is *starlike* if A is compact and there is a point x ∈ A such that for each y ∈ A, the straight-line segment joining x to y lies in A. A subset A of a compact metric space X is *starlike-equivalent* if there is an open neighborhood U of A in X and an embedding e : U → $\mathbb{R}^n$ such that e(U) is an open subset of $\mathbb{R}^n$ and e(A) is a starlike subset of $\mathbb{R}^n$.

**Theorem 2.** Every starlike-equivalent subset of a compact metric space is squeezable.

A proof of this theorem is embedded in the proof of the Lemma in [B]: start at the second paragraph of that proof. The only required change to the argument in [B] is to replace $\mathbb{R}^3$ by $\mathbb{R}^n$, a change which the argument in [B] tolerates without complaint. Alternatively, see the proof of Lemma 5 on pages 55-56 of [D].

## 4. Squashable sets

**Definition.** A subset A of a compact metric space X with metric ρ is a *squashable* subset of X if A is compact and for every null decomposition $\mathcal{G}$ of X such that A ∈ $\mathcal{G}$, every open neighborhood U of A in X, and for every δ > 0, there is a map f ∈ cl($\mathcal{H}_{X-U}(X)$) such that f(A) is a one-point set {p}, f maps X − A homeomorphically onto X − {p}, and for every G ∈ $\mathcal{G}$, either $\text{diam}_\rho(f(G)) < δ$ or f | G = $\text{id}_G$.

The notions of squeezable and squashable are very close. In fact, as the following theorem shows, they are equivalent.

**Theorem 3.** A compact subset of a compact metrizable space X is squeezable if and only if it is squashable.

**Proof.** Let X be a compact metric space with metric ρ, and let A be a compact subset of X.

First assume A is a squeezable subset of X. Let $\mathcal{G}$ be a null decomposition of X such that A ∈ $\mathcal{G}$, let U be an open neighborhood U of A in X, and let δ > 0. We inductively choose open neighborhoods $U_n$ of A in X, real numbers $δ_n$ > 0 and homeomorphisms $h_n$ ∈ $\mathcal{H}_{X-U_n}(X)$ with the following properties.

- Let $U_1$ = U and $δ_1$ = δ/2, and choose $h_1$ ∈ $\mathcal{H}_{X-U_1}(X)$ such that $\text{diam}_\rho(h_1(A)) < δ_1$ and for every G ∈ $\mathcal{G}$, either $\text{diam}_\rho(h_1(G)) < δ_1$ or $h_1$ | G = $\text{id}_G$.

    Let n ≥ 2 and assume that for 1 ≤ i ≤ n − 1, $U_i$ is an open neighborhood of A in X, $δ_i$ > 0 and $h_i$ ∈ $\mathcal{H}_{X-U_i}(X)$ such that $\text{diam}_\rho(h_i(A)) < δ_i$ and for every G ∈ $\mathcal{G}$, either $\text{diam}_\rho(h_i(G)) < δ_i$ or $h_i$ | G = $\text{id}_G$.

- Choose $U_n$ to be an open neighborhood of A in X such that $\text{cl}(U_n) \subset \mathcal{N}_\rho(A, 1/n) \cap U_{n-1}$ and $\text{diam}_\rho(h_{n-1}(U_n)) < δ_{n-1}$.

- Choose $δ_n$ > 0 so that $h_1 \circ h_2 \circ \cdots \circ h_{n-1}$ maps each subset of X of ρ-diameter $< δ_n$ to a subset of X of ρ-diameter $< δ/2^n$.

- Choose $h_n$ ∈ $\mathcal{H}_{X-U_n}(X)$ such that $\text{diam}_\rho(h_n(A)) < δ_n$ and for every G ∈ $\mathcal{G}$, either $\text{diam}_\rho(h_n(G)) < δ_n$ or $h_n$ | G = $\text{id}_G$.

For n ≥ 1, let $k_n = h_1 \circ h_2 \circ \cdots \circ h_n$. Thus, each $k_n$ is a homeomorphism of X.

For i ≥ n ≥ 2, since $h_i \mid X - U_i = id_{X - U_i}$ and $X - U_n \subset X - U_i$, then $k_i \mid X - U_n = k_{n-1} \mid X - U_n$ and $k_i(U_n) = k_{n-1}(U_n)$. Since $diam_\rho(h_{n-1}(U_n)) < \delta_{n-1}$, then $diam_\rho(k_{n-1}(U_n)) < \delta/2^{n-1}$. Hence, $\bar\rho(k_i, k_{n-1}) < \delta/2^{n-1}$ for i ≥ n ≥ 2. Therefore, $\{k_n\}$ is a Cauchy sequence in C(X,X). Hence, $\{k_n\}$ converges to a map $f : X \to X$. Also, the preceding observations in this paragraph imply that $f \mid X - U_n = k_i \mid X - U_n$ and $f(U_n) \subset cl(k_i(U_n)) = k_i(cl(U_n))$ for i ≥ n − 1.

For i ≥ 1, since $X - U \subset X - U_i$, then $h_i \mid X - U = id_{X - U}$. Thus, each $k_n \in \mathcal{H}_{X - U}(X)$. Therefore, $f \in cl(\mathcal{H}_{X - U}(X))$. Note that since a limit of onto maps is onto in C(X,X), then f is onto.

For n ≥ 1, since $f(A) \subset f(U_n) \subset cl(k_{n-1}(U_n))$ and $diam_\rho(cl(k_{n-1}(U_n))) \leq \delta/2^{n-1}$, then $diam_\rho(f(A)) = 0$. Hence, $f(A) = \{p\}$ for some $p \in X$. Since f is onto, then $f(X - A) \supset X - \{p\}$. Since $A = \bigcap_{n \geq 1} cl(U_n)$, then $X - A = \bigcup_{n \geq 1}(X - cl(U_n))$. Hence, $f(X - A) = \bigcup_{n \geq 1} f(X - cl(U_n))$. For each n ≥ 1, since $(X - U_n) \cap cl(U_{n+1}) = \emptyset$, $k_n$ is injective, $f(X - cl(U_n)) \subset k_n(X - U_n)$ and $f(A) \subset k_n(cl(U_{n+1}))$, then $f(X - cl(U_n)) \cap f(A) = \emptyset$. Therefore, $f(X - A) = X - \{p\}$. For each n ≥ 1, $f \mid X - cl(U_n) = k_{n-1} \mid X - cl(U_n)$ is a homeomorphism onto an open subset of $X - \{p\}$. Since $\{ X - cl(U_n) : n \geq 1 \}$ is an increasing sequence of open sets whose union is X − A, it follows that f maps X − A homeomorphically onto X − {p}.

Let $G \in \mathcal{G}$. If G = A, then $diam_\rho(f(G)) = 0$. Suppose G ≠ A. Then there is an n ≥ 1 such that $G \subset X - U_n$. Hence, $f \mid G = k_{n-1} \mid G = h_1 \circ h_2 \circ \cdots \circ h_{n-1} \mid G$. If $h_i \mid G = id_G$ for 1 ≤ i ≤ n − 1, then $f \mid G = id_G$. So assume there is an $i \in \{ 1, 2, \cdots, n - 1 \}$ such that $h_i \mid G \neq id_G$. Let $m = \max \{ i \in \{ 1, 2, \cdots, m \} : h_i \mid G \neq id_G \}$. Then $f \mid G = h_1 \circ h_2 \circ \cdots \circ h_m \mid G$ and $diam_\rho(h_m(G)) < \delta_m$. Therefore, $diam_\rho(f(G)) = diam_\rho(h_1 \circ h_2 \circ \cdots \circ h_{m-1}(h_m(G))) < \delta/2^m < \delta$.

This proves A is a squashable subset of X.

Now assume A is a squashable subset of X. Let $\mathcal{G}$ be a null decomposition of X such that $A \in \mathcal{G}$, let U be an open neighborhood of A in X and let δ > 0. Let V be an open neighborhood of A in X such that $V \subset U$ and such that every element of $\mathcal{G} - \{A\}$ that intersects V has ρ-diameter < δ/3. Since A is squashable, there is a map $f \in cl(\mathcal{H}_{X - V}(X))$ such that $f(A) = \{p\}$ for some $p \in X$, f maps X − A homeomorphically onto X − {p}, and for every $G \in \mathcal{G}$, either $diam_\rho(f(G)) < \delta/3$ or $f \mid G = id_G$. Hence, there is a homeomorphism $h \in \mathcal{H}_{X - V}(X)$ such that $\bar\rho(f,h) < \delta/3$. Since $X - U \subset X - V$, then $h \in \mathcal{H}_{X - U}(X)$. Since $h(A) \subset \mathcal{N}_\rho(\{p\}, \delta/3)$, then $diam_\rho(h(A)) < \delta$. Let $G \in \mathcal{G} - \{A\}$. We must prove that either $diam_\rho(h(G)) < \delta$ or $h \mid G = id_G$. Suppose $h \mid G \neq id_G$. Since $h \mid X - V = id_{X - V}$, then $G \cap V \neq \emptyset$. Hence, $diam_\rho(G) < \delta/3$. Observe that $h(G) \subset \mathcal{N}_\rho(f(G), \delta/3)$. Either $diam_\rho(f(G)) < \delta/3$ or $f \mid G = id_G$. In the first case, $diam_\rho(h(G)) \leq diam_\rho(\mathcal{N}_\rho(f(G), \delta/3)) < \delta$. In the second case, $diam_\rho(h(G)) \leq diam_\rho(\mathcal{N}_\rho(f(G), \delta/3)) = diam_\rho(G, \delta/3)) < \delta$. This proves A is a squeezable subset of X. ∎

We now come to the theorem which justifies the introduction of the notion of squashability: every recursively squashable set is squashable. To deal correctly with the notion of recursive squashability, we need more notation and another definition.

**Notation.** If A is a compact subset of a compact metric space X and $\mathcal{G}_A = \{A\} \cup \{\{x\} : x \in X - A\}$, then $\mathcal{G}_A$ is clearly a decomposition of X. In this situation, we abbreviate $X/\mathcal{G}_A$ to X/A.

**Definition.** If $\mathcal{P}$ is a topological property of compact subsets of compact metrizable spaces (such as squeezable, squashable or starlike-equivalent) and $n \geq 1$, we say that a subset A of a compact metrizable space X is a *recursively $\mathcal{P}$ subset* of X of *filtration length* $\leq n$ if there is a finite sequence $A_1 \subset \cdots \subset A_{n+1}$ of compact subsets of X such that $A_1$ is a $\mathcal{P}$ subset of X, $A_i/A_{i-1}$ is a $\mathcal{P}$ subset of $X/A_{i-1}$ for $2 \leq i \leq n+1$, and $A_{n+1} = A$. Thus, a subset with property $\mathcal{P}$ is recursively $\mathcal{P}$ of filtration length $\leq 0$.

**Theorem 4.** Every recursively squashable subset of a compact metrizable space is squashable.

**Proof.** Clearly, it suffices to prove that if $B \subset A$ are compact subsets of a metric space X such that B is a squashable subset of X and A/B is a squashable subset of X/B, then A is a squashable subset of X. Let $\rho$ be a metric on X, let $\sigma$ be a metric on X/B and let $\pi : X \to X/B$ denote the quotient map. Let $\mathcal{G}$ be a null decomposition of X such that $A \in \mathcal{G}$, let U be an open neighborhood of A in X and let $\delta > 0$. Let $\mathcal{H} = (\mathcal{G} - \{A\}) \cup \{B\} \cup \{\{x\} : x \in B - A\}$. Then $\mathcal{H}$ is a null decomposition of X such that $B \in \mathcal{H}$. Since B is a squashable subset of X, there is a map $g \in \text{cl}(\mathcal{H}_{X-U}(X))$ such that such that $g(B) = \{q\}$ for some $q \in X$, g maps X − B homeomorphically onto X − {q}, and for every $G \in \mathcal{G}$, either $\text{diam}_\rho(g(G)) < \delta$ or $g \mid G = \text{id}_G$.

Observe that $h = \pi \circ g^{-1} : X \to X/B$ is a homeomorphism. $\pi(\mathcal{G}) = \{\pi(G) : G \in \mathcal{G}\}$ is a null decomposition of X/B such that $\pi(A) = A/B \in \pi(\mathcal{G})$. Also $\pi(U)$ is an open neighborhood of A/B in X/B. There is an $\epsilon > 0$ such that $h^{-1}$ maps each subset of X/B of $\sigma$-diameter < $\epsilon$ to a subset of X of $\rho$-diameter < $\delta$. Since A/B is a squashable subset of X/B, there is a map $f \in \text{cl}(\mathcal{H}_{X/B - \pi(U)}(X/B))$ such that $f(A/B) = \{p\}$ for some $p \in X/B$, f maps X/B − A/B homeomorphically onto X/B − {p}, and for every $G \in \mathcal{G}$, either $\text{diam}_\sigma(f(\pi(G))) < \epsilon$ or $f \mid \pi(G) = \text{id}_{\pi(G)}$.

Define the map $k : X \to X$ by $k = h^{-1} \circ f \circ h \circ g$. We now run through the details of the verification that k has the correct properties to imply that A is squashable in X. First observe that since $h \circ g = \pi$, then $k = h^{-1} \circ f \circ \pi$.

**Lemma.** If X, Y and Z are compact metrizable spaces, then composition $(f,g) \mapsto g \circ f : C(X,Y) \times C(Y,Z) \to C(X,Z)$ is continuous.

**Proof.** Let σ and τ be metrics on Y and Z, respectively. Let $\bar{\sigma}$, $\bar{\tau}$ and $\bar{\bar{\tau}}$ denote the supremum metrics on C(X,Y), C(Y,Z) and C(X,Z) determined by σ, τ and τ, respectively. Let f ∈ C(X,Y), g ∈ C(Y,Z) and ϵ > 0. There is a δ > 0 such that δ < ϵ/2 and: y, y′ ∈ Y and σ(y,y′) < δ ⇒ τ(g(y),g(y′)) < ϵ/2. Now suppose f′ ∈ C(X,Y) and g′ ∈ C(Y,Z) such that $\bar{\sigma}$(f,f′) < δ and $\bar{\tau}$(g,g′) < δ. Since $\bar{\sigma}$(f,f′) < δ, then $\bar{\bar{\tau}}$(g∘f,g∘f′) < ϵ/2; and $\bar{\bar{\tau}}$(g∘f′,g′∘f′) ≤ $\bar{\tau}$(g,g′) < δ < ϵ/2. Hence, $\bar{\bar{\tau}}$(g∘f,g′∘f′) ≤ $\bar{\bar{\tau}}$(g∘f,g∘f′) + $\bar{\bar{\tau}}$(g∘f′,g′∘f′) < ϵ. This proves (f,g) ↦ g∘f is continuous. ∎

Since g ∈ cl($\mathcal{H}_{X-U}$(X)) and f ∈ cl($\mathcal{H}_{X/B-\pi(U)}$(X)), there are sequences {$g_n$} in $\mathcal{H}_{X-U}$(X) and {$f_n$} in $\mathcal{H}_{X/B-\pi(U)}$(X/B) such that {$g_n$} converges to g in C(X,X) and {$f_n$} converges to f in C(X/B,X/B). Define the homeomorphism $k_n$ : X → X by $k_n$ = $h^{-1}$∘$f_n$∘h∘$g_n$. Since composition is continuous by the preceding lemma, then {$k_n$} converges to k. Since $g_n$ | X − U = $id_{X-U}$, h(X − U) = X/B − π(U) and $f_n$ | X/B − π(U) = $id_{X/B-\pi(U)}$, then $k_n$ | X − U = $id_{X-U}$. Thus, {$k_n$} ⊂ $\mathcal{H}_{X-U}$(X) and k ∈ cl($\mathcal{H}_{X-U}$(X)).

k(A) = $h^{-1}$∘f∘π(A) = $h^{-1}$∘f(A/B) = $h^{-1}$({p}) which is a one-point set in X. Since π maps X − A homeomorphically onto X/B − A/B, f maps X/B − A/B homeomorphically onto X/B − {p}, and $h^{-1}$ : X/B → X is a homeomorphism, then k maps X − A homeomorphically onto X − $h^{-1}$(p).

$diam_\rho$(k(A)) = $diam_\rho$($h^{-1}$({p})) = 0. Let G ∈ $\mathcal{G}$ − {A}. Then G ∈ $\mathcal{H}$. Then k(G) = $h^{-1}$∘f∘π(G). Either $diam_\sigma$(f(π(G))) < ϵ or f | π(G) = $id_{\pi(G)}$. If $diam_\sigma$(f(π(G))) < ϵ, then $diam_\rho$(k(G)) = $diam_\rho$($h^{-1}$(f∘π(G))) < δ. Suppose f | π(G) = $id_{\pi(G)}$. Then k | G = $h^{-1}$∘f∘π | G = $h^{-1}$∘π | G = g | G. Since G ∈ $\mathcal{H}$, then either $diam_\rho$(g(G)) < δ or g | G = $id_G$. Thus, either $diam_\rho$(k(G)) = $diam_\rho$(g(G)) < δ or k | G = $id_G$. ∎

### 5. Consequences

Theorems 1, 3 and 4 imply:

**Corollary 1.** If $\mathcal{G}$ is a null decomposition of a compact metric space X each of whose non-singleton elements is recursively squeezable, then $\mathcal{G}$ is strongly shrinkable (and, thus, the quotient map π : X → X/$\mathcal{G}$ is a strong near homeomorphism).

Then using Theorem 2, we have:

**Corollary 2.** If $\mathcal{G}$ is a null decomposition of a compact metric space X each of whose non-singleton elements is recursively starlike-equivalent, then $\mathcal{G}$ is strongly shrinkable (and, thus, the quotient map π : X → X/$\mathcal{G}$ is a strong near homeomorphism).

Clearly Corollary 2 implies the following result which is the main theorem of [MOR].

**Theorem 1.1 of [MOR].** If n ≥ 1 and $\mathcal{G}$ is a null decomposition of a compact metric space X each of whose non-singleton elements is recursively starlike-equivalent subset of X of filtration length ≤ n, then $\mathcal{G}$ is strongly shrinkable (and the quotient map π : X → X/$\mathcal{G}$ is a strong near homeomorphism).